\documentclass[12pt]{article}

\usepackage{amsfonts,amsmath,amssymb,amsthm}
\def\lb{\label}
\def\s#1{T_{#1}}
\def\l{t}

\def\gins#1{\rule{-7pt}{0pt} {\small #1} \rule{-10.5pt}{0pt} }
\def\qTr#1{{\rm Tr_{#1}}}
\def\oTr#1{{\cal T}\kern-.1667em r^{#1}}
\def\qdim{{\rm qdim}}
\def\qed{\hfill \rule{5pt}{5pt}}

\newcommand{\non}{\nonumber}
\newcommand{\ts}{\,}
\newcommand{\tss}{\hspace{1pt}}

\newcommand{\CC}{\mathbb{C}}
\newcommand{\Sym}{\mathfrak S}
\newcommand{\Tc}{{\mathcal T}}
\newcommand{\Uc}{{\tss\mathcal U}}
\newcommand{\Vc}{{\mathcal V}}
\newcommand{\He}{\mathcal H}
\newcommand{\qin}{q^{-1}}
\newcommand{\al}{\alpha}

\newcommand{\de}{\delta}

\newcommand{\la}{\lambda}

\newcommand{\si}{\sigma}

\newcommand{\beq}{\begin{equation}}
\newcommand{\eeq}{\end{equation}}
\newcommand{\ben}{\begin{equation*}}
\newcommand{\een}{\end{equation*}}

\newcommand{\bpf}{\begin{proof}}
\newcommand{\epf}{\end{proof}}

\def\beql#1{\begin{equation}\label{#1}}

\newtheorem{thm}{Theorem}[section]
\newtheorem{lemma}[thm]{Lemma}
\newtheorem{prop}[thm]{Proposition}
\newtheorem{cor}[thm]{Corollary}

\theoremstyle{definition}
\newtheorem{definition}[thm]{Definition}
\newtheorem{example}[thm]{Example}

\theoremstyle{remark}
\newtheorem{remark}[thm]{Remark}

\newcommand{\bth}{\begin{thm}}
\renewcommand{\eth}{\end{thm}}
\newcommand{\bpr}{\begin{prop}}
\newcommand{\epr}{\end{prop}}
\newcommand{\ble}{\begin{lemma}}
\newcommand{\ele}{\end{lemma}}
\newcommand{\bco}{\begin{cor}}
\newcommand{\eco}{\end{cor}}
\newcommand{\bex}{\begin{example}}
\newcommand{\eex}{\end{example}}
\newcommand{\bde}{\begin{definition}}
\newcommand{\ede}{\end{definition}}
\newcommand{\bre}{\begin{remark}}
\newcommand{\ere}{\end{remark}}

\begin{document}
\title{On the idempotents of Hecke algebras}
\author{A.P. Isaev${\tss}^{a}$,  A.I. Molev${\tss}^{b}$
and A.F. Os'kin${\tss}^{a}$}
\date{\empty}
\maketitle

\begin{center}
${}^{a}$ Bogoliubov Laboratory of Theoretical Physics,
Joint Institute for Nuclear Research, \\
Dubna, Moscow region 141980, Russia \\
E-mail: isaevap@theor.jinr.ru,  oskinandrej@gmail.com
\\
\vspace{0.3cm}
${}^{b}$ School of Mathematics and Statistics  \\
University of Sydney, NSW 2006, Australia\\
E-mail: alexm@maths.usyd.edu.au
\end{center}

\vspace{0.3cm}

\begin{abstract}
We give a new construction of primitive idempotents
of the Hecke algebras associated with the symmetric groups.
The idempotents are found as evaluated products
of certain rational functions thus providing a new version
of the fusion procedure for the Hecke algebras.
We show that the normalization
factors which occur in the procedure
are related to the Ocneanu--Markov trace of the idempotents.
\end{abstract}

\section{Introduction}

It was observed by Jucys~\cite{j:yo} that the primitive idempotents 
of the symmetric group $\Sym_n$ can be obtained by taking certain 
limit values of the rational function \beql{phiu} 
\Phi(u_1,\dots,u_n)= \prod_{1\leqslant i<j\leqslant n} 
\Big(1-\frac{(i\ts j)}{u_i-u_j}\Big), \eeq where $u_1,\dots,u_n$ are 
complex variables and the product is calculated in the group algebra 
$\CC[\Sym_n]$ in the lexicographical order on the pairs $(i,j)$. A 
similar construction, now commonly referred to as the {\it fusion 
procedure\/}, was developed by Cherednik~\cite{c:sb}, while complete 
proofs were given by Nazarov~\cite{n:yc}. A simple version of the 
fusion procedure establishing its equivalence with the Jucys--Murphy 
construction was recently found by one of us in \cite{Mol}; see also 
\cite[Ch.~6]{m:yc} for applications to the Yangian representation 
theory and more references. In more detail, let $\Tc$ be a standard 
tableau associated with a partition $\lambda$ of $n$ and let 
$c_k=j-i$, if the element $k$ occupies the cell of the tableau in 
row $i$ and column $j$. Then the consecutive evaluations \ben 
\Phi(u_1,\dots,u_n)\big|_{u_1=c_1}\big|_{u_2=c_2}\dots 
\big|_{u_n=c_n} \een are well-defined and this value yields the 
corresponding primitive idempotent $E^{\lambda}_\Tc$ multiplied by 
the product of the hooks of the diagram of $\lambda$. The left ideal 
$\CC[\Sym_n]\ts E^{\lambda}_\Tc$ is the irreducible representation 
of $\Sym_n$ associated with $\lambda$, and the $\Sym_n$-module 
$\CC[\Sym_n]$ is the direct sum of the left ideals over all 
partitions $\lambda$ and all $\lambda$-tableaux $\Tc$.

Our aim in this paper is to derive an analogous version
of the fusion procedure for the Hecke algebra $\He_n=\He_n(q)$
associated with $\Sym_n$. The procedure goes back to
Cherednik~\cite{Cher}, while detailed proofs
relying on $q$-versions of the Young symmetrizers
were given by Nazarov~\cite{nazar};
see also Grime~\cite{g:hfq} for its hook
version. We use a different approach based on
the formulas for the primitive idempotents of $\He_n$
in terms of the {\it Jucys--Murphy elements\/}.
These formulas derived by Dipper and James~\cite{dj:bi}
generalize the results
of Jucys~\cite{j:fy} and Murphy~\cite{mu:is} for $\Sym_n$.

The main result of this paper is an explicit formula
for the orthogonal primitive idempotents of $\He_n$.
The idempotents are obtained as
a result of consecutive evaluations of a rational
function similar to \eqref{phiu}.
The normalization factors in the expressions
for the Hecke algebra idempotents turn out to be
related to the Ocneanu--Markov trace
of the idempotents.

\section{Idempotents of $\He_n$ and Jucys--Murphy elements}

Let $q$ be a formal variable.
The Hecke algebra $\He_{n}$ over the field $\CC(q)$
is generated by the elements $T_1,\dots,T_{n-1}$
subject to the defining relations
\begin{align}
\non
T_i \, T_{i+1} \, T_i &=
T_{i+1} \, T_i \,  T_{i+1},\\
T_{i} \,  T_{j} &= T_{j} \,  T_{i}\quad
{\rm for} \quad |i-j| > 1,
\non\\
T^2_i &= 1 + (q-\qin)\ts  T_i.
\non
\end{align}
Given a reduced decomposition $w=s_{i_1}\dots s_{i_l}$
of an element $w\in\Sym_n$ in terms of the generators
$s_i=(i,i+1)$, set $T_{w}=T_{i_1}\dots T_{i_l}$. Then
$T_{w}$ does not depend on the reduced decomposition,
and the set $\{T_w\ |\ w\in\Sym_n\}$ is a basis of
$\He_n$ over $\CC(q)$.

The {\it Jucys-Murphy elements\/}
$y_1,\dots,y_n$ of $\He_n$ are defined inductively by
\begin{gather}
y_1 = 1,\qquad
y_{k+1} = \s{k}\ts y_{k}\ts \s{k}
\qquad\text{for}\quad k=1,\dots,n-1.\label{ym}
\end{gather}
These elements satisfy
\ben
y_k\ts\s{m}=\s{m}\ts y_k,\qquad m \neq k,k-1.
\een
In particular, $y_1,\dots,y_n$ generate
a commutative subalgebra of $\He_n$.
The elements $y_k$ can be written
in terms of the elements $T_{(i\ts j)}\in\He_n$,
associated with the transpositions $(i\ts j)\in\Sym_n$ as follows:
\ben
y_k=1+(q-\qin)\ts
\big(T_{(1\ts k)}+T_{(2\ts k)}+\dots+T_{(k-1\ts k)}\big).
\een
Hence, the normalized elements $(y_k -1)/(q-\qin)$
specialize to the Jucys--Murphy elements for $\Sym_n$
as $q\to 1$; see \cite{j:fy}, \cite{mu:is}, \cite{dj:bi}.

For any $k=1,\dots,n$ we let $w_k$ denote the
unique longest element
of the symmetric group $\Sym_k$ which is regarded as
the natural subgroup of $\Sym_n$.
The corresponding elements $T_{w_k}\in\He_n$
are then given by $T_{w_1}=1$ and
\begin{align}
T_{w_k} &= \s1(\s2\, \s1) \cdots
(\s{k-2}\dots \s1)(\s{k-1}\s{k-2}\dots \s{1})
\lb{gars} \\[0.2cm]
{}&=    (\s1 \dots \s{k-2}\, \s{k-1})
(\s1 \dots \s{k-2}) \cdots (\s1\, \s2)\ts\s1,\qquad
k=2,\dots,n.
\lb{gars1}
\end{align}
We point out the following properties of the
elements $T_{w_k}$ which are easily verified by induction
with the use of \eqref{gars} and \eqref{gars1}:
\begin{align}
T_{w_k}\ts \s{j} &= \s{k-j}\ts T_{w_k}, \qquad\quad
1\leqslant j<k\leqslant n,
\label{gf3} \\
T_{w_k}^2 &= y_1 y_2 \cdots y_{k},\qquad  k=1,\dots,n.
\non
\end{align}
Following \cite{nazar}, for each $i=1,\dots,n-1$ set
\begin{gather}
T_i(x, y) = \dfrac{T_i \, y - T^{-1}_i \, x}{y - x}
=T_i+\frac{q-\qin}{x^{-1}y-1},
\label{bxt}
\end{gather}
where $x$ and $y$ are complex variables. We will regard
the $T_i(x, y)$ as rational functions in $x$ and $y$ with values in
$\He_n$. It is well-known that
they satisfy the relations
\begin{equation}
\lb{ybeH}
T_i(x,y) \, T_{i+1}(x,z) \, T_i(y,z) =
T_{i+1}(y,z) \, T_i(x,z) \, T_{i+1}(x,y),
\end{equation}
(the Yang--Baxter equation), and
\begin{gather}
\label{baxtinverse}
T_i(x,y)\ts T_i(y,x)=
\dfrac{\left(x - q^2 y \right)\left(x - q^{-2} y\right)}
{\left( x-y \right)^2}.
\end{gather}

\ble\label{lem:tom}
We have the identities
\beql{prop2}
T_{w_{k}} \, \s{j}(x, y) = \s{k-j}(x, y) \, T_{w_{k}},
\qquad 1 \leqslant j <k\leqslant n,
\eeq
and
\beql{prop4}
T_{w_{k+1}}\, \s2(u, \si_{k-1}) \dots \s{k}(u,\si_1)\,
T_{w_{k}}^{-1} =
 T_{w_{k}}\, \s1(u,\si_{k-1}) \dots \s{k-1}(u,\si_1)\,
 T_{w_{k-1}}^{-1}
\s{k},
\eeq
where $1\leqslant k< n$ and
$u,\si_1,\dots,\si_{k-1}$ are complex parameters.
\ele

\bpf
Relation (\ref{prop2}) is immediate from (\ref{gf3}), while
(\ref{prop4}) is deduced from
\ben
 (\s{k} \dots \s2\s1)\s{j}(x,y)  =
\s{j-1}(x,y) \, (\s{k} \dots\s2 \s1), \qquad 2
\leqslant j \leqslant k,
\een
by taking into account
the identity
\ben
T_{w_{k}}^{-1} T^{}_{w_{k+1}} = T_{w_{k-1}}^{-1} T^{}_k T^{}_{w_{k}} =
\s{k} \dots \s2 \s1
\een
implied by (\ref{gars}) and (\ref{gars1}).
\epf

Now we recall the construction of the orthogonal primitive
idempotents for the Hecke algebra from \cite{dj:bi}.
We will identify
a partition $\lambda=(\lambda_1,\dots,\lambda_l)$ of $n$
with its diagram which is
a left-justified array of rows of cells such that the top row
contains $\lambda_1$ cells, the next row contains
$\lambda_2$ cells, etc.
A cell outside $\la$ is called {\it addable\/} to $\lambda$
if the union of $\lambda$ and the cell is a diagram.
A tableau $\Tc$ of shape $\lambda$ (or a $\lambda$-tableau $\Tc$)
is obtained by
filling in the cells of the diagram
bijectively with the numbers $1,\dots,n$.
A tableau $\Tc$ is called standard if its entries
increase along the rows and down the columns.
If a cell occurs in row $i$ and column $j$,
its $q$-{\it content\/} will be defined as $q^{2\tss(j-i)}$.

In accordance to \cite{dj:bi}, a set of orthogonal
primitive idempotents $\{E^{\la}_{\Tc}\}$ of $\He_n$, parameterized
by partitions $\la$ of $n$ and standard $\la$-tableaux $\Tc$
can be constructed inductively by the following rule.
Set $E^{\la}_\Tc=1$ if $n=1$, whereas for $n\geqslant 2$,
\beql{murphyfo}
E^{\la}_\Tc=E^{\mu}_\Uc\ts
\frac{(y_n-\rho_1)\dots (y_n-\rho_k)}{(\si-\rho_1)\dots (\si-\rho_k)},
\eeq
where $\Uc$ is the tableau
obtained from $\Tc$ by removing the
cell $\al$ occupied by $n$, $\mu$ is the shape of $\Uc$, and
$\rho_1,\dots,\rho_k$ are the $q$-contents of all addable cells of $\mu$
except for $\alpha$, while $\si$ is the $q$-content of the latter.
In particular, if $\la$ and $\la'$ are partitions of $n$, then
\ben
E^{\la}_\Tc\ts E^{\la'}_{\Tc'}
=\de_{\la\la'}\de_{\Tc\Tc'}\ts E^{\la}_\Tc
\een
for arbitrary standard tableaux $\Tc$ and $\Tc'$ of shapes
$\la$ and $\la'$, respectively.
Moreover,
\ben
\sum_{\la}\sum_\Tc E^{\la}_\Tc=1,
\een
summed over all partitions $\la$ of $n$ and all standard
$\la$-tableaux $\Tc$.

In what follows we will omit the superscript $\la$ and write
simply $E_\Tc$ instead of $E^{\la}_\Tc$.
Given a standard
$\la$-tableau $\Tc$ and $k\in\{1,\dots,n\}$, we set
$\si_k=q^{2\tss(j-i)}$ if the element $k$ of $\Tc$ occupies
the cell in row $i$ and column $j$. Then
\begin{align}\label{jmev}
y_k \, E_\Tc =  E_\Tc \, y_k = \si_k\ts E_\Tc.
\end{align}
Furthermore, given a standard
tableau $\Uc$ with $n-1$ cells, the corresponding idempotent
$E_\Uc$ can be written as
\begin{align}
E_\Uc = \sum_{\Tc} E_{\Tc},
\label{brr}
\end{align}
summed over all standard tableaux $\Tc$ obtained from $\Uc$ by
adding one cell with entry $n$.
Exactly as in the case of the symmetric group $\Sym_n$ (see \cite{Mol}),
this relation can be used to derive the following alternative form
of \eqref{murphyfo}. Consider the rational function
\begin{align}\label{fubrr}
E_\Tc(u)=E_\Uc\ts \dfrac{u - \si_n}{u - y_n}
\end{align}
in a complex variable $u$ with values in $\He_n$.
Then this function is regular at $u=\si_n$ and the corresponding
value coincides with $E_\Tc$:
\begin{align}\label{idem1}
E_\Tc = E_\Uc\ts \dfrac{u - \si_n}{u - y_n}\Big|_{u=\si_n}.
\end{align}

\section{Fusion formulas for primitive idempotents}

For $k=1,\dots, n-1$ introduce the elements of $\He_n$ by
\begin{align}
\label{yyyy}
Y_{k}(\si_1, \si_2, \dots, \si_k; u)
= T_{w_{k}}\, \s{k}(\si_1, u) \s{k-1}(\si_2, u) \dots
\s1(\si_k, u)\, T_{w_{k+1}}^{-1},
\end{align}
where $\si_1, \si_2, \dots, \si_k$ and $u$ are complex parameters.

\ble\label{clem}
Let $\Uc$ be a standard tableau with $k$ cells and the $q$-contents
$\si_1,\si_2, \dots, \si_k$. Then
\begin{align}
\label{prop5}
&E_\Uc\ts Y_k(\si_1, \dots, \si_k; u) =  \nonumber\\
&=
 (u- \si_1) \left(\prod_{j = 1}^k \frac
{\left(u - q^2 \si_j \right) \left(u - q^{-2} \si_j\right)}
{\left(u - \si_j\right)^2}
\right) E_\Uc(u - y_{k+1})^{-1}.
\end{align}
\ele

\bpf
We start with representing (\ref{prop5}) in the form
\begin{gather}
\label{prr}
\left(u-\si_1\right)^{-1} E_\Uc (u - y_{k+1}) =
E_\Uc \, T_{w_{k+1}} \, \s1(u,\si_k) \dots
\s{k}(u,\si_1) \, T_{w_{k}}^{-1},
\end{gather}
where we have used~\eqref{baxtinverse} and taken into account
the fact that $E_\Uc$ commutes with $y_{k+1}$.
Now we prove (\ref{prr}) by induction. For $k=1$ we have
\ben
(u-\si_1)^{-1}(u - \s1^{2}) = T_{1} \cdot \s1(u,\si_1),
\een
which is true, as $\si_1=1$.
Due to (\ref{prop2}) and (\ref{prop4}), the right hand side of
(\ref{prr}) can be written in the form
\begin{multline}
E_\Uc\, \s{k}(u,\si_k)\, T_{w_{k+1}}\,
\s2(u,\si_{k-1})\dots\s{k}(u,\si_1)\,
T_{w_{k}}^{-1} = \\
= E_\Uc\, \s{k}(u,\si_k)\, T_{w_{k}}\, \s1(u,\si_{k-1}) \dots
\s{k-1}(u,\si_1)\, T_{w_{k-1}}^{-1} \s{k}.
\non
\end{multline}
Using \eqref{brr}, we can
write $E_\Uc = E_\Uc E_\Vc$, where $\Vc$ is
the tableau obtained from $\Uc$ by removing the cell occupied by
$k$. Hence, the right hand side of
(\ref{prr}) becomes
\begin{multline}
E_\Uc E_\Vc \, \s{k}(u,\si_k)\,T_{w_{k}}\, \s1(u,\si_{k-1})
 \dots \s{k-1}(u,\si_1)\, T_{w_{k-1}}^{-1} \s{k} = \\
 = E_\Uc \, \s{k}(u,\si_k)\, \Big(E_\Vc T_{w_{k}}\, \s1(u,\si_{k-1})
\dots \s{k-1}(u,\si_1)\, T_{w_{k-1}}^{-1}\Big) \s{k} = \\
= \left(u - \si_1\right)^{-1} E_\Uc \s{k}(u,\si_k)(u - y_k) \s{k}.
\non
\end{multline}
The last equality holds by the induction hypothesis. Now we represent
$\s{k}(u,\si_k)$ in the form
\ben
\s{k}(u,\si_k) = \dfrac{\s{k} \, \si_k -
\s{k}^{-1}\, u}{\si_k - u} =  \s{k} +
\frac{(q-\qin)\ts u}{\si_k - u}.
\een
This gives
\begin{multline}
E_\Uc \s{k}(u,\si_k)(u - y_k)\ts \s{k}
= E_\Uc \left(\s{k} + \frac{(q-\qin)\ts u}{\si_k - u} \right)
\left(u - y_k\right)\s{k} = \\
= E_\Uc \left( -u\ts (q-\qin)\ts \s{k} + u\ts \s{k}^{2} - y_{k+1}\right)
= E_\Uc \left(u - y_{k+1}\right),
\non
\end{multline}
thus completing the proof.
\epf

Let $\la=(\la_1,\dots,\la_l)$ be a partition of $n$.
We will use the conjugate partition $\la'=(\la'_1,\dots,\la'_m)$
so that $\la'_j$ is the number of cells in the $j$-th column of $\la$.
If $\al=(i,j)$ is a cell of $\la$, then the corresponding hook
is defined as $h_{\al}=\la_i+\la'_j-i-j+1$ and the content
is $c_{\al}=j-i$. Set
\beql{star3}
f(\la)=\prod_{\al\in\la}\frac{\ q^{\ts c_{\al}}\ }{[h_{\al}]_q},
\eeq
where we have used the notation
\ben
[n]_q = \frac{q^n -q^{-n}}{q-q^{-1}}.
\een

Suppose that
$\Tc$ is a standard $\la$-tableau.
As before, for each $k\in\{1,\dots,n\}$
we let $\si_k$ denote the $q$-content $q^{2(j-i)}$ of the cell
$(i,j)$ occupied by $k$ in $\Tc$. Consider the rational function
\ben
F_{n} (u) = \dfrac{u - \si_{n}}{u - \si_1}\ts \prod_{k = 1}^{n - 1}
\dfrac{(u - \si_k)^2}{(u - q^2 \si_k)(u - q^{-2}\si_k)}.
\een

\ble\label{lem:norm}
The rational function $F_{n}(u)$ is regular at $u=\si_n$
and
\ben
F_{n}(\si_n)=f(\mu)^{-1}\ts f(\la),
\een
where $\mu$ denotes the shape of the
standard tableau obtained from $\Tc$
by removing the cell occupied by $n$.
\ele

\bpf
It is clear that $F_{n}(u)$ depends only on the shape $\mu$
and does not depend on the standard tableau $\Uc$ obtained from
$\Tc$ by removing the cell occupied by $n$. Therefore,
we may assume that $\Uc$ is the row tableau obtained by writing
the elements $1,\dots,n-1$ into the cells of $\mu$ consecutively
by rows starting with the top row.
Suppose that
the rows of $\mu$ are
\ben
\mu_1=\dots=\mu_{p_1}>\mu_{p_1+1}=\dots=\mu_{p_2}>\dots>
\mu_{p_{s-1}+1}=\dots=\mu_{p_s}
\een
for some integers $p_1,\dots,p_s$
such that $1\leqslant p_1<p_2<\dots<p_s$
and some $s\geqslant 1$.
With this notation, $F_{n} (u)$ can be written in
the form
\ben
F_{n} (u) = (u - \si_{n})\ts \prod_{i = 1}^{s}
(u - q^{2\ts\mu_{p_i}-2\ts p_i})\ts \prod_{i = 0}^{s}
(u - q^{2\ts \mu_{p_i+1}-2\ts p_i})^{-1},
\een
where we set $p_0=0$ and $\mu_{p_s+1}=0$. Possible values
of the $q$-content $\si_n$ are
$\si_n=q^{2\ts\mu_{p_j+1}-2\ts p_j}$ for $j=0,1,\dots,s$.
Hence, for a fixed value of $j$ the factor $u - \si_{n}$
cancels, and so $F_{n} (\si_n)$
is well-defined
and can be expressed in the form
\beql{fnsi}
F_{n} (\si_n)=
\big(\ts q^{2\ts\mu_{p_j+1}}-q^{2\ts\mu_{p_j+1}+2}\ts\big)\ts
\prod_{\al\in\mu}(1-q^{2h_{\al}})
\prod_{\al\in\la}(1-q^{2h_{\al}})^{-1},
\eeq
which is verified by a simple calculation.
On the other hand, $f(\la)$ can be represented as
\ben
f(\la)=q^{b(\la)}\ts(1-q^2)^n
\prod_{\al\in\la}(1-q^{2h_{\al}})^{-1},\qquad
b(\la)=\sum_{i\geqslant 1}\la_i(\la_i-1).
\een
Therefore, the expression in
\eqref{fnsi} equals $f(\mu)^{-1}\ts f(\la)$,
as required.
\epf

Introduce the rational function $\Psi(u_1,\dots,u_n)$
in complex variables $u_1,\dots,u_n$ with values in $\He_n$
by the formula
\ben
\Psi(u_1,\dots,u_n)=\prod_{k=1,\dots,\tss n-1}^{\longrightarrow}
\Big(\s{k}(u_1, u_{k+1})\ts \s{k-1}(u_2, u_{k+1}) \dots
\s1(u_k, u_{k+1})\Big)\cdot T_{w_{n}}^{-1}.
\een
As before, we let
$\la$ be a partition of $n$ and let
$\Tc$ be a standard $\la$-tableau.

\bth\label{theo}
The idempotent $E_{\Tc}$ can be obtained by
the consecutive evaluations
\beql{theor}
E_{\Tc} = f(\la)\cdot
\Psi(u_1,\dots,u_n)
\big|_{u_1=\si_1} \big|_{u_2=\si_2} \dots \big|_{u_n=\si_n},
\eeq
where the rational functions are regular at
the evaluation points at each step.
\eth

\bpf
We argue by induction on $n$. For $n\geqslant 2$ we
let $\Uc$ denote the standard tableau obtained from $\Tc$
by removing the cell occupied by $n$ and let $\mu$ be the shape
of $\Uc$.
Applying Lemma~\ref{lem:norm} and
the induction hypothesis, we can write the right hand side of
\eqref{theor} in the form
\ben
F_{n}(\si_n) \,
E_{\Uc} \, Y_{n-1}(\si_1, \dots, \si_{n-1}; u_{n})\big|_{u_n=\si_n},
\een
where the elements $Y_{n-1}(\si_1, \dots, \si_{n-1}; u_{n})$
are defined in \eqref{yyyy}. The proof is completed
by the application of Lemma~\ref{clem} and
relation~\eqref{idem1}.
\epf

\bex
Using~\eqref{theor}, for $n = 3$ and $\la=(2,1)$ we get
\begin{align}
\label{ia}
E_\Tc = \frac{1}{[3]_q} \s1(\si_1, \si_2) \s2(\si_1, \si_3)
\s1(\si_2, \si_3) (\s1\s2\s1)^{-1}.
\end{align}
In particular,
\ben
\si_1=1,\quad  \si_2=q^2,\quad  \si_3=q^{-2}\qquad \text{for}\qquad
\Tc= \begin{tabular}{|c|c|}
\hline
\gins{1} & \gins{2} \\
\hline
\gins{3 }& \multicolumn{1}{c}{}\\
\cline{1-1}
\end{tabular}
\een
and
\ben
\si_1=1,\quad  \si_2=q^{-2},\quad  \si_3=q^2\qquad \text{for}\qquad
\Tc= \begin{tabular}{|c|c|}
\hline
\gins{1} &\gins{3} \\
\hline
\gins{2 }& \multicolumn{1}{c}{}\\
\cline{1-1}
\end{tabular}\ .
\een
Note that (\ref{ia}) can be reduced to the fusion formulas contained
in \cite[p.~106]{Isa1}.
\qed
\eex

\bex
For $n = 4$ and $\la=(2^2)$ the idempotent $E_\Tc$ is obtained
by evaluating the rational function
\begin{gather}
\label{t4}
\dfrac{1}{[3]_q [2]_q^2}\ts
\s1(u_1, u_2)
\s2(u_1, u_3) \s1(u_2, u_3)
\s3(u_1, u_4)\s2(u_2, u_4)\s1(u_3, u_4)\ts T_{w_4}^{-1}
\end{gather}
consecutively at $u_1=\si_1$, $u_2=\si_2$, $u_3=\si_3$, and $u_4=\si_4$.
We have
\ben
\si_1 = 1,\quad \si_2 = q^2,\quad \si_3 = q^{-2}, \quad \si_4 = 1\qquad
\text{for}
\qquad \Tc =
\begin{tabular}{|c|c|}
\hline
\gins{1} & \gins{2} \\
\hline
\gins{3} & \gins{4} \\
\hline
\end{tabular}
\een
and
\ben
\si_1 = 1,\quad \si_2 = q^{-2},\quad \si_3 = q^{2}, \quad \si_4 = 1\qquad
\text{for}
\qquad \Tc =
\begin{tabular}{|c|c|}
\hline
\gins{1} & \gins{3} \\
\hline
\gins{2} & \gins{4} \\
\hline
\end{tabular}\ \   .
\een
Note that for both tableaux the expression (\ref{t4})
contains the factor $\s3(u_1, u_4)$
which is not defined for $u_1=\si_1$ and $u_4=\si_4$. Nevertheless,
the whole expression (\ref{t4}) is regular under the
consecutive evaluations due to Theorem~\ref{theo}.
We will use this example to illustrate
the relationship with the approach of \cite{nazar}.
Using the
relation (\ref{ybeH}) one can rewrite (\ref{t4}) as
\ben
\dfrac{1}{[3]_q [2]_q^2}\ts
\s2(u_2, u_3)  \s1(u_1, u_3) \s2(u_1, u_2)
\s3(u_1, u_4)\s2(u_2, u_4)\s1(u_3, u_4)\ts T_{w_4}^{-1}.
\een
By \cite[Lemma~2.1]{nazar},
the product $\s2(u_1, u_2) \s3(u_1, u_4) \s2(u_2, u_4)$ is equal to
\begin{multline}
\frac{\left( (T_2 u_2 - T_2^{-1} u_1)\ts T_3\ts (T_2 u_4 - T_2^{-1} u_2)
+(q-\qin)\ts u_1 ((q-\qin)\ts u_2 T_2  + u_2-u_1) \right)}
{(u_2-u_1)(u_4-u_2)}\\
{}-\frac{(q-\qin)\ts u_1 (u_{1} - q^2 u_{2})
(u_1 - q^{-2} u_2)}{(u_2-u_1)(u_4-u_1)(u_4-u_2)}
\nonumber
\end{multline}
and it is regular
for $u_1 = q^{\pm 2}u_2$ at $u_1 = u_4$. It was shown
in \cite{nazar} that such
considerations can be extended to the general
expression (\ref{theor}) to prove that it is
regular in the limits
$u_i \to \si_i$.
\qed
\eex

We conclude this section by showing that taking an appropriate
limit in Theorem~\ref{theo}
as $q\to 1$ we can recover the respective formulas of
\cite{Mol} for the primitive idempotents
of the symmetric group $\Sym_n$.

Take the parameters $x$ and $y$ in \eqref{bxt} in the form
$x = q^{2\tss u}$ and $y = q^{2\tss v}$.
Since $\s{i} \xrightarrow[q \to 1]{} s_{i}$,
for the limit
value of  $\s{i}(x, y)$ we have
\begin{gather}\label{bxtyang}
\s{i}(x, y) = \s{i} + \frac{q^{u - v}}{[v - u]_q}
\xrightarrow[q \to 1]{} s_{i}\, \varphi_{i, i+1} (u, v),
\end{gather}
where
\ben
\varphi_{i,j}(u, v) = 1 - \frac{(i\ts j)}{u - v}.
\een
Using (\ref{bxtyang}) we can calculate the
corresponding limit for the element (\ref{yyyy}) to get
\begin{gather}
\label{yyyy2}
Y_{k}(\si_1, \si_2, \dots, \si_k; u)  \xrightarrow[q \to 1]{}
\varphi_{1, k+1} (c_1, u) \varphi_{2, k+1} (c_2, u) \dots
\varphi_{k, k+1} (c_k, u),
\end{gather}
where $\si_m=q^{2\tss c_m}$.
Clearly, the normalization factor $f(\la)$ specializes
to the inverse of the product of the hooks of $\la$, and so
the substitution of (\ref{yyyy2})
into (\ref{theor}) leads to the main result of \cite{Mol}.

\section{The Ocneanu--Markov trace of the idempotents}

The purpose of this section is to calculate the
Ocneanu--Markov trace of the idempotents $E_{\Tc}$
which turns out to be related to the
normalization factor $f(\la)$ defined in (\ref{star3}).

\bde
For any given standard tableau $\Tc$ with $n$
cells, its {\it quantum dimension\/} is defined as
\begin{gather}
\label{qdim5}
\qdim\ts \Tc = \oTr{n}(E_\Tc) \; ,
\end{gather}
where $\oTr{n}: \He_n \to\CC$ is
the Ocneanu--Markov trace; see e.g. \cite{Jon1}.
\qed
\ede

The Ocneanu--Markov trace $\oTr{n}$
can be defined as
the composition of the maps
\ben
\oTr{n} =
\qTr{1}\qTr{2}\dots\qTr{n}.
\een
The linear maps
$\qTr{m+1}: \He_{m+1}\to \He_{m}$ from the Hecke
algebra $\He_{m+1}$ to its natural subalgebra $\He_{m}$
are determined by the following properties,
where $Q\in \CC$ is a fixed parameter, while
$X, Y \in \He_m$ and $Z \in \He_{m+1}$:
\begin{align}
\label{qtrdef}
& \qTr{m+1}(XZY) = X\qTr{m+1}(Z)Y, \quad \qTr{m+1}(X) = Q X, \nonumber\\
&\qTr{m+1}(\s{m}^{\pm 1} X \s{m}^{\mp1})
= \qTr{m}(X),\quad \qTr{m+1}(\s{m}) =
1,\\
&\qTr{m}\qTr{m+1}(\s{m}Z) = \qTr{m}\qTr{m+1}(Z\s{m}).\nonumber
\end{align}

Our calculation of (\ref{qdim5}) is based on the approach
of~\cite{IsOgH}.
The following statement can be found in that paper.

\bpr\label{zprop}
Consider the rational function
in $u$ with values in the Hecke algebra $\He_m$
which is defined by
\ben
Z_{m+1}(u) = \qTr{m+1}\left(u - y_{m+1}\right)^{-1},\qquad
y_{m+1} \in \He_{m+1},
\een
where $\He_m$ is regarded as a subalgebra of $\He_{m+1}$.
Then,
\begin{gather}\label{zpropeq}
Z_{m+1}(u) = \dfrac{\l\ts Q + u-1}{t\ts u(u - 1)}
\left( \prod_{k = 1}^{m}
\dfrac{(u - y_k)^2}{(u - q^2 y_k)(u - q^{-2} y_k)}
- \dfrac{(1-\l\ts Q)(u-1)}{\l\ts Q + u-1} \right),
\end{gather}
where $\l=q-\qin$.
\epr

\bpf
>From the definition of the Jucys-Murphy
elements (\ref{ym}) we deduce the identity
\begin{gather}
\label{zmeq2}
\dfrac{1}{u - y_{m+1}} = \s{m}\dfrac{1}{u - y_{m}}\s{m}^{-1}
+ \dfrac{1}{u - y_m} \left( \s{m}^{-1}
+ \dfrac{\l u}{(u - y_{m+1})} \right)
\dfrac{\l y_m}{(u - y_m)}.
\end{gather}
Applying the map Tr$_{m+1}$ to both sides of~\eqref{zmeq2}
and using \eqref{qtrdef}
we get
\ben
\dfrac{(u - q^2 y_m)(u - q^{-2} y_m)}{(u - y_m)^2}
Z_{m+1}(u) = Z_m(u) + \dfrac{\l
(1 - Q \l) y_m}{(u - y_m)^2}.
\een
For all $k =1,\dots,m+1$ introduce the function $\bar{Z}_k(u)$ by
\ben
Z_{k}(u) = \bar{Z}_k(u) + (Q -
\l^{-1})u^{-1}.
\een
This gives the relation
\ben
\bar{Z}_{m+1}(u) = \dfrac{(u - y_m)^2}{(u - q^2 y_m)
(u - q^{-2} y_m)}\bar{Z}_m(u).
\een
Solving this recurrence relation with the initial condition
\ben
\bar{Z}_1(u) = \qTr1(u - y_1)^{-1} - (Q - \l^{-1})u^{-1} =
\dfrac{t\ts Q + u - 1}{t\ts u(u - 1)},
\een
we come to \eqref{zpropeq}.
\epf

The normalization factor $f(\la)$ defined in
\eqref{star3} and
the quantum dimension~(\ref{qdim5})
turn out to be related as shown in the following proposition.
As before, we let $\la$ be
a partition of $n$, and $\Tc$ a standard $\la$-tableau.

\bpr\label{prop:flaeq}
We have the relation
\ben
f(\la) = \qdim\ts \Tc\ts \prod_{k = 1}^{n}\si_k
\left(Q + \frac{\si_k - 1}{q-\qin}\right)^{-1}.
\een
\epr

\bpf
Using~\eqref{fubrr} and~\eqref{idem1} we get
\ben
\qTr{n} (E_{\Tc}) = \qTr{n}E_\Tc(u)\big|_{u=\si_n} = E_\Uc
\ts (u - \si_{n})\qTr{n}(u - y_{n})^{-1}\big|_{u=\si_n}.
\een
Using equations~\eqref{zpropeq} and taking into account~\eqref{jmev}
we obtain
\begin{multline}
\non
\qTr{n} (E_\Tc) = \frac{1}{\si_n}\left(Q + \frac{\si_{n} - 1}{\l}\right)
E_\Uc\\
{}\times\dfrac{u -
\si_{n}}{u - 1}\left(\prod_{k = 1}^{n-1}
\dfrac{(u - \si_k)^2}{(u - q^2 \si_k)(u - q^{-2} \si_k)} -
(u - 1)\dfrac{1 - \l Q}{\l Q+u-1}\right)\Bigg|_{u=\si_n} = \\
= \frac{1}{\si_n}\left( Q + \frac{\si_{n} - 1}{\l}\right) E_\Uc
F_{n}(\si_n).
\end{multline}
Applying the maps $\qTr{k}$ consequently,
we finally obtain
\ben
\qdim\ts \Tc =
\oTr{n}(E_\Tc) = \qTr1\qTr2\dots\qTr{n} (E_\Tc) =
\prod_{m=1}^{n}\frac{1}{\si_m}\left( Q + \frac{\si_{m} -
1}{\l}\right)\ts F_{m} (\si_m).
\een
The statement now follows from Lemma~\ref{lem:norm}.
\epf

The following corollary is immediate from
Proposition~\ref{prop:flaeq}.

\bco
The Ocneanu--Markov trace $\oTr{n}(E_{\Tc})$
of the idempotent $E_{\Tc}$ depends only on the shape
$\la$ of $\Tc$ and does not depend on $\Tc$.
\eco

\vspace{0.2cm}
\noindent
{\bf Acknowledgements}

The work of the first author was partially supported
by the RFBR grant No. 08-01-00392-a. The second author
gratefully acknowledges the support
of the Australian Research Council.

\end{document}